\documentclass[reqno]{amsart}
\usepackage{etoolbox}
\usepackage[utf8]{inputenc}
\usepackage[margin=24mm,headheight=15pt]{geometry}
\usepackage{amsmath,amssymb,amsthm,mathtools,bm}
\usepackage[]{amsmath}
\usepackage[numbers,sort&compress]{natbib}
\usepackage[hidelinks]{hyperref}
\usepackage[]{amssymb}
\usepackage{graphicx}
\usepackage{float}
\usepackage{expl3}
\usepackage{indentfirst}
\usepackage{ mathrsfs }
\usepackage{mathtools}
\usepackage{cmap}
\usepackage{exercise}
\usepackage{tikz}
\usetikzlibrary{tikzmark,calc}
\usetikzlibrary{decorations.pathmorphing}
\usepackage{breqn}
\usepackage{cancel}
\usepackage{comment}
\usepackage{xcolor}
\usepackage{tcolorbox}

\newtheorem{theo}{Theorem}[section]

\newtheorem{ex}[theo]{Example}
\newtheorem{rem}{Remark}
\newtheorem{defin}[theo]{Definition}
\newtheorem{lemm}[theo]{Lemma}
\newtheorem{coro}[theo]{Corollary}
\newtheorem{prop}[theo]{Proposition}
\numberwithin{equation}{section}

\newcommand{\D}{\mathcal D}
\newcommand{\I}{\mathcal I}
\newcommand{\OO}{\mathcal O}
\newcommand{\x}{\mathbf x}
\newcommand{\y}{\mathbf y}
\newcommand{\Res}{\operatorname{Res}}
\newcommand{\ee}{\mathbf e}
\newcommand{\W}{\mathcal W}
\newcommand{\C}{\mathbb C}
\newcommand{\sgn}{\operatorname{sgn}}

\title{Differential equations for Macdonald's hypergeometric series via Dunkl operators of type A}
\author{Fan Liu}
\address{Beijing International Mathematics Research Center, Peking University, Beijing 100871, China}
\email{liufan\_math@pku.edu.cn}

\begin{document}

\begin{abstract}
We solve an eigenvalue problem for a family of trigonometric Dunkl operators depending on a parameter, extending the result of Ref.\cite{Dunkl25} to the case of a general parameter. We introduce families of nonsymmetric $W$-type operators closely related to Dunkl operators and analyze their algebraic structure, various representations and symmetrization. We then use these operators to characterize non-symmetric Macdonald's hypergeometric series introduced in Ref.\cite{Brennecken22} and recover several known results. These results suggest a direct connection between Macdonald's hypergeometric series and Dunkl operators.
\end{abstract}
\maketitle
\tableofcontents

\section{Introduction}
The generalized classical hypergeometric series is defined by \cite{Erdelyi53} 
\begin{align}\label{rfs}
	_rF_s(a_1,\cdots,a_r;b_1,\cdots,b_s;x)&=\sum_{k=0}^{\infty}\frac{\prod_{i=1}^r(a_i)_k}{\prod_{j=1}^s(b_j)_k}\frac{x^k}{k!},\quad r,s\in \mathbb{N}.
\end{align} 
Especially, $_2F_1(a,b;c;x)$ is the well-known Gauss hypergeometric function. At a formal level, the generalized hypergeometric series can be obtained from the power-series expansion $e^x=\sum_{n=0}^{\infty}\frac{x^n}{n!}$ by inserting suitable Pochhammer factors $(a)_{n}=\prod_{k=0}^{n-1}(a+k)$ into the place of numerator or denominator. Such a series can be characterized by the following linear ordinary differential equation (ODE). 
\begin{align}\label{ODE}
	\left(x\frac{d}{dx}\prod_{j=1}^{s} (x\frac{d}{dx}+b_j-1)-\prod_{i=1}^{r}(x\frac{d}{dx}+a_i-1)x\right)F(x)&=0.
\end{align}
In general, the series (\ref{rfs}) is not the unique solution of the ODE. A complete description of the ODE’s solution space requires a careful analysis of its monodromy representation.

There are many ways \cite{Heckman87,Gelfand90,Lauricella} to generalize hypergeometric functions to several variables. Herz \cite{Herz55} introduced hypergeometric series in matrix arguments, and Constantine \cite{Constantine63} subsequently expanded them in zonal polynomials. Later, Macdonald \cite{Mac13I,Mac13II} introduced two families of hypergeometric series expanded respectively in Jack and Macdonald polynomials in his manuscripts. This generalization is natural, since Jack and Macdonald polynomials \cite{Stanley89,Mac95} can be viewed as one- and two-parameter deformations of zonal spherical polynomials. We refer to these families as Macdonald’s $\beta$-deformed and $(q,t)$-deformed hypergeometric series, respectively, although similar series were introduced and studied by several other authors around the same time \cite{Koranyi91,Yan92,Kaneko93,Kaneko96}.

The Gauss-type Macdonald's $\beta$-deformed hypergeometric series has received the most extensive study. Yan \cite{Yan92} and Kaneko \cite{Kaneko93} independently derived constraint equations for this series and obtained its Selberg integral representation by establishing existence and uniqueness of solutions to those equations. Kaneko also derived the holomonic differential system for the Gauss-type series and later extended these methods to the $(q,t)$-deformed case \cite{Kaneko96}.

Around the same time, during his study of spherical harmonics for measures invariant under finite reflection groups, Dunkl \cite{Dunkl89} introduced a family of commuting differential–difference operators. These Dunkl operators were later extended to the trigonometric version in various ways \cite{Heckman91,Cherednik91,Polychronakos92} and ultimately unified within the framework of Cherednik algebras \cite{Cherednik05}. In general, a commuting family of operators admits joint eigenfunctions. This principle led to extensive theories of orthogonal polynomials associated with root systems, including Heckman’s hypergeometric functions \cite{Heckman87} and the Koornwinder polynomials \cite{Koornwinder,Mac03}. All of them can be regarded as the  multivariable generalizations of hypergeometric orthogonal polynomials in the Askey–Wilson scheme \cite{Koekoek}. From a physical perspective, suitably modified Dunkl operators can be used to producing higher Hamiltonians for certain quantum Calogero–Moser–Sutherland many-body systems \cite{Calogero69,Sutherland71,Olshanetsky83} and their various integrable deformations \cite{Ruij,Diejen,Bernard93}.

The theories of Macdonald's hypergeometric series and Dunkl operators can be brought together naturally. By studying the algebra of Dunkl operators of type $A_{N-1}$, Baker and Forrester \cite{Baker97,Baker00} constructed hypergeometric generating functions for multivariable Hermite, Laguerre and Al-Salam–Carlitz polynomials. They used these generating functions to solve a number of combinatorial problems. In their work \cite{Baker98}, the Dunkl's integral kernel \cite{Dunkl91} arises naturally as a non-symmetric extension of Macdonald’s $\beta$-deformed hypergeometric series. More recently, Brennecken and R\"osler \cite{Brennecken22} have provided a complete definition for non-symmetric hypergeometric series and carried out an in-depth study on Dunkl-Laplace transformation. These non-symmetric series may be a more fundamental object since they can be symmetrized to recover the original symmetric series.

Determining constraint equations and integral representations for Macdonald's hypergeometric series with general parameters has long been an important problem, and recent work has led to substantial progress. To obtain operator representations of the partition functions of certain $\beta$- and $(q,t)$-deformed matrix models \cite{Wang22}, we introduced families of homogeneous $W$-type difference and differential operators \cite{Liu24}. These families arise from polynomial representations of the quantum toroidal algebra \cite{Feigin12} and spherical degenerate double-affine Hecke algebra (dDAHA) \cite{Schiffmann13}, respectively. More recently, these $W$-operators have been used to derive constraint equations for Macdonald's hypergeometric series \cite{Liu25,Liu26}. Chen and Sahi \cite{Chen25} established conditions ensuring the uniqueness of solutions to these equations.
Using these uniqueness theorems, we derived the integral representations for certain hypergeometric series. These results suggest a direct connection between the Macdonald's hypergeometric series and the superintegrability of matrix models \cite{Mironov22,Cassia20,Cassia25}.

Several questions nevertheless remain open. First, the explicit forms of these $W$-operators are quite complicated and require simplification. Mironov and Morozov \cite{Mironov23,Mironov24} have expressed the $W$-operators in terms of Dunkl operators. We will provide a rigorous proof of this representation here. Second, the constraint equations for non-symmetric Macdonald hypergeometric series remain an unexplored problem. Finally, we would like to derive a holomonic differential system for Macdonald hypergeometric series, as an extension of the Kaneko's results on Gauss-type case. In this paper, we focus on Macdonald’s $\beta$-deformed hypergeometric series. We introduce families of non-symmetric $W$-operators by using Dunkl operators and study their properties in order to address the problems outlined above.

The paper is organized as follows. In Section 2, we review symmetric and non-symmetric Jack polynomials and Dunkl operators. In Section 3, we introduce a family of Heckman–Polychronakos–Dunkl operators and obtain their eigenvalues. In Section 4, we introduce families of non-symmetric $W$-operators and analyzes their algebraic structure and various representations. In Section 5, we use these $W$-operators to give a characterization of non-symmetric Macdonald's $\beta$-deformed hypergeometric series and recover some previous known results.

\section{Preliminaries}
Throughout this paper, we fix a positive integer $N$ and let $\beta$ and $\gamma$ are two real parameters. We mainly refer to Ref.\cite{Mac95} for the part of partitions and symmetric functions and Ref.\cite{Baker98} for the part of Dunkl operators and non-symmetric Jack polynomials.
\subsection{Partial order and partitions}
We denote the space $V_N=\mathbb{N}^N$ generated by $N$-tuple $\eta=(\eta_1,\cdots,\eta_N)$ of
non-negative integers. The weight $|\eta|$ and length $\ell(\eta)$ are defined by
\begin{align}
	|\eta|=\sum_{i=1}^{N}\eta_i,\qquad \ell(\eta)=\sharp\{i:\eta_i>0\}.
\end{align}
A partition refers to a tuple $\lambda=(\lambda_1,\lambda_2,\ldots)$ of
non-negative integers in weakly decreasing order. Then we denote the set of partitions with length less than $N$ as
\begin{align}
	\mathcal{P}_N:=\{(\lambda_1,\ldots,\lambda_N)\in V_N:\lambda_1\ge\cdots \ge\lambda_N\}.
\end{align}
For a given $\eta\in V_N$, we denote $\eta_+$ is the unique element in $\{\mathfrak{S}_N\cdot \eta\}\cap \mathcal{P}_N$.

The partial order ``$\prec$" on $V_N$ is defined by
\begin{align}
	\zeta\prec \eta \ \text{iff}\ \begin{cases}
	\zeta^+<_D \eta^+ & \zeta^+\neq \eta^+,\\
	\zeta<_D \eta & \zeta^+= \eta^+,
	\end{cases}
\end{align}
where ``$\le_D$" is the usual dominance order on $V_N$, i.e.,
\begin{align}\label{order}
	\zeta\le_D \eta \ \text{iff} \ |\zeta|=|\eta|  \text{and}\ \sum_{i=1}^r \zeta_i\le\sum_{i=1}^r \eta_i,\ \forall r=1,\cdots,N,
\end{align}
and $\zeta<_D\eta$ means $\zeta\le_D \eta$ and $\zeta\neq\eta$. 

By adding and removing a basic unit in the $i$-th row from $\lambda\in\mathcal{P}_N$ respectively, one will get the two $N$-tuples as
\begin{align*}
	\lambda_{(i)}&=(\lambda_1\cdots,\lambda_{i-1},\lambda_{i}-1,\cdots,\lambda_{i+1},\cdots,\lambda_N), \quad \lambda_i>0,\nonumber\\
	\lambda^{(i)}&=(\lambda_1\cdots,\lambda_{i-1},\lambda_{i}+1,\cdots,\lambda_{i+1},\cdots,\lambda_N),
\end{align*}
which are not belong to $\mathcal{P}_N$ anymore in general.

For a given $\lambda\in \mathcal{P}_N$, one can identify its Young diagram as $\{(i,j)\in\mathbb{Z}^2: 1\le i\le N, 1\le j\le \lambda_i\}$. For $a\in\C$, the $\beta$-Pochhammer symbol is defined by
\begin{align}\label{contants}
	[a]_{\lambda}&:=\prod_{(i,j)\in\lambda}[a+(j-1)+(1-i)\beta]=\prod_{i=1}^N(a+(1-i)\beta)_{\lambda_i}.
\end{align}

\subsection{Symmetric polynomials}
Let $\x=(x_1,\cdots,x_N)$ and $\y=(y_1,\cdots,y_N)$. Then we denote $\C[\x]=\C[x_1,\ldots,x_N]$ the polynomial ring in $N$ variables and $\C[\x]^{\mathfrak{S}_N}$ subring of symmetric polynomials, where $\mathfrak{S}_N$ is the permutation group on $N$ letters.

For $\lambda\in\mathcal{P}_N$, the monomial symmetric polynomials $m_{\lambda}$, the power-sums $p_\lambda$ and elementary symmetric polynomials $e_{\lambda}$ are respectively defined by 
\begin{align}
	m_{\lambda}(\x)&=\sum_{\eta\in \mathfrak{S}_N\cdot \lambda}x^\eta,\quad  x^\eta=\prod_{i=1}^Nx_i^{\eta_i},\\
	p_{\lambda}(\x)&=\prod_{i=1}^Np_{\lambda_i}(\x), \quad  p_n(\x)=\sum_{i=1}^Nx_i^n,\\
	e_{\lambda}(\x)&=\prod_{i=1}^Ne_{\lambda_i}(\x),\quad  e_n(\x)=\sum_{1\le i_1\le \cdots\le i_n \le N}\prod_{k=1}^Nx_{i_1}\cdots x_{i_n}.
\end{align}
All of them are the basis of $\C[\x]^{\mathfrak{S}_N}$. Here $\{e_n\}_{n=0,\cdots,N}$ are usually characterized by its generating function 
\begin{align}
	\sum_{n=1}^Ne_n(\x)z^{N-n}=\prod_{i=1}^N(x_i+z).
\end{align}

Let us consider the hermitian inner product on the polynomial space $\mathbb{C}[\mathbf{x}]$
\begin{align}\label{inner}
	(f,g):=\int_{\mathbb{T}^N} f(\mathbf{x})g(\mathbf{x}^{-1})\prod_{i\neq j}|1-x_i/x_j|^{\beta} d\x,
\end{align}
where $\mathbb{T}^{N}=\{(x_1,\cdots,x_N)\in\mathbb{C}^N;|x_i|=1,\ i=1,\cdots,N\}$ and $d\mathbf{x}=\frac{dx_1}{x_1}\wedge \cdots \wedge \frac{dx_N}{x_N}$.
The Jack polynomial $P_{\lambda}$ are uniquely determined by the following conditions
\begin{subequations}\label{cond}
\begin{align}
	& P_{\lambda}=m_{\lambda}+\sum_{\mu<_D\lambda}u_{\lambda,\mu}m_{\mu}, \label{cond1} \\
	& (P_{\lambda}, m_{\mu})=0, \quad \forall \mu\in \mathcal{P}_N, \ \mu<_D\lambda, \label{cond2}
\end{align}
\end{subequations}
where $u_{\lambda,\mu}$ are $\beta$-deformed analogues of the Kostka coefficients. 

For our purposes, we need to consider a different normalization of Jack polynomials as 
\begin{align}
	C_{\lambda}=\frac{|\lambda|!}{\prod_{(i,j)\in\lambda}[(\lambda_i-j+1)+\beta(\lambda_j'-i)]}P_{\lambda}.
\end{align}
such that $\sum_{\lambda\in\mathcal{P}_N, |\lambda|=n}C_{\lambda}(\x)=(p_1(\x))^n$ and 
\begin{align}
	C_{\lambda}(1^N)=\frac{|\lambda|![N\beta]_{\lambda}}{\prod_{(i,j)\in\lambda}(c_i-j+\beta\lambda_j')(c_i-j+\beta\lambda_j'+1-\beta)}.
\end{align}

The degree-one Pieri formulas for Jack polynomials $C_{\lambda}$ are given by \cite{Lassalle90}
\begin{subequations}\label{Pieri}
\begin{align}
	p_1(\x)C_{\lambda}(\x) &= \frac{1}{1+|\lambda|}\sum_{i=1,\lambda^{(i)}\in\mathcal{P}_N}^{\ell(\lambda)+1} \binom{\lambda^{(i)}}{\lambda}_{\beta}C_{\lambda^{(i)}}(\x), \\
	\sum_{i=1}^N\frac{\partial}{\partial x_i}\cdot \frac{C_{\lambda}(\x) }{C_{\lambda}(1^N) }&= \sum_{i=1,\lambda_{(i)}\in\mathcal{P}_N}^{\ell(\lambda)} \binom{\lambda}{\lambda_{(i)}}_{\beta}\frac{C_{\lambda_{(i)}}(\x) }{C_{\lambda_{(i)}}(1^N) }
\end{align}
\end{subequations}
for any $\ell(\lambda)\le N$, where the generalized binomial coefficients are given by 
\begin{align}\label{binom}
\binom{\lambda}{\lambda_{(i)}}_{\beta}
=[c_i+\beta(\ell(\lambda)-1)]
\prod_{\substack{j=1\\ j\ne i}}^{\ell(\lambda)}
\frac{c_j - c_i - \beta}{c_j - c_i},
\end{align}
and 
$c_i=c_i(\lambda)=\begin{cases}
	\lambda_i-(i-1)\beta & i\le \ell(\lambda),\\
	-(i-1)\beta & i>\ell(\lambda).
\end{cases}$

\subsection{Dunkl operators and non-symmetric Jack polynomials}
The Dunkl operators of type $A_{N-1}$ are given by
\begin{align}\label{Dunkl}
	\D_i=\D_i(\x):=\partial_i+\beta\sum_{j=1,j\neq i}^N\frac{1}{x_i-x_j}(1-s_{ij}),\quad i=1,\cdots,N,
\end{align}
where $\partial_i=\frac{\partial}{\partial x_i}$ and $s_{ij}\in\mathfrak{S}_N$ is the operator permuting the variables $x_i$ and $x_j$. These Dunkl operators are commutative and covariant, i.e., 
\begin{align}
	[\D_i,\D_j]&=0, \\
	w\D_iw^{-1}&=\D_{w(i)},
\end{align}
for $i,j=1,\cdots,N$ and $w\in\mathfrak{S}_N$.

We define an inner product on $\C[\x]$ by
\begin{align}\label{noninner}
	(f,g)':=\rho(f)g(\x)|_{\x=0^N},
\end{align} 
where $\rho$ is the isomorphism defined by $\rho(f(\x))=f(\D_1,\cdots,\D_N)$ and $\x=a^N$ means $x_i=a$ for $i=1,\cdots,N$ and $a\in\C$. For a given operator $X$, we denote $X^{\dagger}$ is the adjoint with respected to $(\cdot,\cdot)'$, i.e., $(Xf,g)'=(f,X^{\dagger}g)'$.
Then, it is clear that
\begin{align}\label{adjoint2}
	(x_i)^{\dagger}=\D_i, \ (\D_i)^{\dagger}=x_i, \quad i=1,\cdots,N.
\end{align}

The trigonometric analogous of the Dunkl operators are various. We denote
\begin{subequations}\label{CHPdunkl}
\begin{align}
	\pi_i&:=x_i\D_i,\\
	\pi^{H}_i&:=x_i\D_i+\beta\sum_{j=1,j\neq i}^N(1-s_{ij}),\\
	d_i&:=x_i\D_i-\beta\sum_{j=1,j>i}^N(1-s_{ij}),
\end{align}
\end{subequations}
for $i=1,\cdots,N$. The operators $\pi_i$, $\pi_i^{H}$ and $d_i$ are respectively introduced by Polychronakos \cite{Polychronakos92}, Heckman \cite{Heckman91} and Cherednik \cite{Cherednik91}. 
\begin{lemm}\label{selfadjoint}\cite{Heckman91,Dunkl91}
(1) All of the three sets of operators in (\ref{CHPdunkl}) are self-adjoint with respect to the inner products (\ref{inner}) and (\ref{noninner}).\\
(2) The operators $\{d_i\}_{i=1,\cdots,N}$ are commutative but not covariant, while $\{\pi_i,\pi_i^{H}\}_{i=1,\cdots,N}$ are covariant but not commutative, i.e., for $i,j=1,\cdots,N$ and $w\in\mathfrak{S}_N$, 
\begin{subequations}
\begin{align}
	[d_i,d_j]&=0,\label{CDop} \\
	w\pi_iw^{-1}&=\pi_{w(i)},\\
	w\pi_i^{H}w^{-1}&=\pi^{H}_{w(i)}.
\end{align}
\end{subequations}
\end{lemm}
The Cherednik-Dunkl operators $\{d_i\}_{i=1,\cdots,N}$ generate a commutative subalgebra of the dDAHA of $GL(N)$. Bernard \cite{Bernard93} used the polynomial operators $\{e_n(d_1,\cdots,d_N)\}_{n=1,\cdots,N}$ to describe the conservation quantities of spin trigonometric Calogero-Sutherland model.

The non-symmetric Jack polynomials $E_{\eta}$ for $\eta\in V_N$ are uniquely determined by the conditions
\begin{subequations}
\begin{align}
	& E_{\eta}=x^{\zeta}+\sum_{\zeta\prec\eta}u_{\eta,\zeta}x^{\zeta},  \\
	& (E_{\eta}, x^{\zeta})=0, \quad \forall \zeta\in V_N, \ \zeta\prec\eta.
\end{align}
\end{subequations}
Here one can get also same polynomials by replacing the second condition with 
the orthogonality $(E_{\eta}, x^{\zeta})'=0$ for all $\zeta\prec\eta$ or the eigen equation 
\begin{align}\label{dinonJ}
	d_iE_{\eta}(\x)=[q_i(\eta)-(1-i)]\beta E_{\eta}(\x), \quad i=1,\cdots,N,
\end{align} 
where $q_i(\eta)=\beta^{-1}\eta_i-\sharp\{j>i:\eta_j>\eta_i\}-\sharp\{j<i:\eta_j\ge \eta_i\}$.
The action of the permutation operator $s_i=s_{i,i+1}$ on the non-symmetric Jack polynomial is given by \cite{Baker98}
\begin{align}
s_i E_\eta =
\begin{cases}
	\displaystyle
	\frac{1}{q_i(\eta)-q_{i+1}(\eta)} E_\eta
	+\left(1-\frac{1}{(q_i(\eta)-q_{i+1}(\eta))^2}\right) E_{s_i\eta},
	& \eta_i > \eta_{i+1}, \\
	E_\eta,
	& \eta_i = \eta_{i+1}, \\
	\displaystyle
	\frac{1}{q_i(\eta)-q_{i+1}(\eta)} E_\eta + E_{s_i\eta},
	& \eta_i < \eta_{i+1}.
\end{cases}
\end{align}
It means that
\begin{align}\label{SnonJ}
	wE_{\eta}\in\left\{\sum_{\zeta}f_{\zeta}E_{\zeta}:f_{\zeta}=0 \ \text{unless}\ \zeta^+=\eta^+\right\}.
\end{align}
for any $w\in\mathfrak{S}_N$ and $\eta\in V_N$.

We sometimes apply another normalization 
\begin{align}\label{nonL}
	L_{\eta}=\frac{|\eta|!}{\prod_{(i,j)\in\eta}(\eta_i-j+1+\beta\l_{\eta}(i,j))}E_{\eta},
\end{align}
where $\l_{\eta}(i,j)=\sharp\{k>i:j\le\eta_k\le\eta_i\}+\sharp\{k<i:j\le\eta_k+1\le\eta_i\}$. The polynomial $L_{\eta}$ satisfies
$\sum_{\eta\in V_N,|\eta|=n}L_{\eta}=(p_1(\x))^n$. Then the Jack polynomials and the non-symmetric Jack polynomials are connected by
\begin{align}\label{CLR}
	C_{\lambda}=\sum_{\eta\in V_{N},\eta^+=\lambda}L_{\eta}.
\end{align}

\section{Eigenvalues of Heckman-Polychronkos-Dunkl operators}
\subsection{Heckman-Polychronkos-Dunkl operators}
Due to (\ref{dinonJ}), the eignvalues of the Cherednik-Dunkl operators $\{d_i\}_{i=1,\cdots,N}$ are known. We intend to study the Dunkl operators introduced by Heckman and Polychronkos and their eigenvalue problem. Hence, we define
\begin{align}\label{HPDunkl}
	\pi_i(\gamma):=\pi_i(\gamma;\x)=x_i\partial_i
	+\frac\beta2\sum_{1\le j\ne i\le N}\left(\frac{x_i+x_j}{x_i-x_j}+\gamma\right)(1-s_{ij}),\quad i=1,\cdots,N,
\end{align}
where $\partial_i=\frac{\partial}{\partial x_i}$ and $s_{ij}\in \mathfrak{S}_N$ is the operator permuting $x_i$ and $x_j$. The operator $\pi_i(\gamma)$ reduces to $\pi^{H}$ and $\pi$ in \eqref{CHPdunkl} when $\gamma=0$ and $\gamma=1$, respectively. 
\begin{defin}
For $n\in \mathbb{N}$, we define 
\begin{equation}\label{Inc}
	\I_n(\gamma):=\Res\sum_{i=1}^N\pi_i^n(\gamma),
\end{equation} 
where $\Res(f)$ means the action of $f$ is restraint on the symmetric polynomial ring $\C[\x]^{\mathfrak{S}_N}$.
\end{defin}
These operators depends on an additional parameter $\gamma$ and recover the Heckman–Polychronakos operators of Ref.\cite{Dunkl25} at $\gamma=1$. We therefore refer to $\{{\I_n(\gamma)}_{n\in\mathbb{N}}\}$ as the Heckman–Polychronakos–Dunkl operators. 
\begin{ex}
We list the operator $\I_n(\gamma)$ for $n=1$ and $2$
\begin{align}
	\I_1(\gamma)&=\sum_{i=1}^Nx_i\partial_i,\\
	\I_2(\gamma)&=\sum_{i=1}^{N}\left(x_i\partial_i\right)^2 + \beta \sum_{i \neq j} \frac{x_i+x_j}{x_i-x_j}\left(x_i\partial_i-x_j\partial_j\right).
\end{align}
\end{ex}
$\I_1(\gamma)$ and $\I_2(\gamma)$ are respectively the momentum and (gauged) Hamiltonian of the quantum Sutherland model \cite{Sutherland71}, respectively, and both of them are independent to $\gamma$. 

Now, we will prove the commutativity of $\I_n(\gamma)$ by following Heckman's method in Ref.\cite{Heckman91}.

\begin{prop}
For a given $\lambda\in \mathcal{P}_{N}$, the linear polynomial space
\begin{align}\label{Gamma}
	\Gamma(\lambda):=\{f=\sum_{\eta\in V_N}f_{\eta}x^{\eta}: f_{\eta}=0 \ \text{unless} \ \eta\in C(\lambda)\}
\end{align}
is invariant under the action of $\pi_i(\gamma)$ for $i=1,\cdots,N$, where the Weyl chamber $C(\lambda)$ is given by
\begin{align}
	C(\lambda)=\{\eta\in V_N: w(\eta)\le_D \lambda, \forall w\in \mathfrak{S}_N\}.
\end{align}
\end{prop}

\begin{proof}
For $p,q\in\mathbb{N}$, we denote $\min=\min\{p,q\}$ and $\max=\max\{p,q\}$. Then we have
\begin{align}\label{action}
	&\left(\frac{x_i+x_j}{x_i-x_j}+\gamma\right)(1-s_{ij})x_i^px_j^q\nonumber\\
	=&\sgn(p-q)\left(\sum_{r=\min}^{\max-1}+\sum_{r=\min+1}^{\max}\right)x_i^{r}x_j^{p+q-r}+\gamma(x_i^{p}x_j^{q}-x_i^{q}x_j^{p}),
\end{align}
where $\sgn(x)=\begin{cases}
	1 & x> 0\\
	0 & x= 0\\
	-1 & x<0\\
\end{cases}$ is the sign function. The result means $\pi_i(\gamma)x^\eta$ belongs to $\Gamma(\lambda)$ in (\ref{Gamma}) for $\eta\in C(\lambda)$ and $i=1,\cdots,N$. This completes the Proposition.
\end{proof}

Since the partition $\lambda\in \mathcal{P}_N$ is highest element in the orbit $\mathfrak{S}_N\cdot \lambda$ in terms of the dominance order (\ref{order}),  we arrive at
\begin{coro}\label{invariant}
For $\lambda\in \mathcal{P}_{N}$ fixed, the linear space
\begin{align}
	\Gamma_+(\lambda):=\{f=\sum_{\mu\in \mathcal{P}_{N}}f_{\mu}m_{\mu}(\x): f_{\mu}=0\ \text{unless}\ \mu\le_D\lambda\}
\end{align}
is invariant under the Heckman–Polychronakos–Dunkl operators $\{\I_{n}(\gamma)\}_{n\in\mathbb{N}}$.
\end{coro}

\begin{theo}
The Heckman–Polychronakos–Dunkl operators $\{\I_{n}(\gamma)\}_{n\in\mathbb{N}}$ act diagonally on the Jack polynomials.
\end{theo}
\begin{proof}
According to (\ref{cond1}), we get
\begin{align}
	\Gamma_+(\lambda):=\{f=\sum_{\mu\in \mathcal{P}_{N}}f'_{\mu}P_{\mu}(\x); f'_{\mu}=0\ \text{unless}\ \mu\le_D\lambda\}.
\end{align}
Then follows from Corollary \ref{invariant}, we have that
$\I_{n}(\gamma)$ act triangularly on $m_{\lambda}(\x)$ and $P_{\lambda}(\x)$, i.e.,
\begin{subequations}\label{Inm}
\begin{align}
	\I_n(\gamma)m_{\lambda}(\x)&=\sum_{\mu\in\mathcal{P}_N,\mu\le_D\lambda}b^{(n)}_{\lambda,\mu}m_{\mu}(\x), \\
	\I_n(\gamma)P_{\lambda}(\x)&=\sum_{\mu\in \mathcal{P}_N, \mu\le_D\lambda}c^{(n)}_{\lambda,\mu}P_{\nu}(\x), 
\end{align}
\end{subequations}
where $b^{(n)}_{\lambda,\mu}$ and $c^{(n)}_{\lambda,\mu}$ are certain coefficient matrices related to $\beta$ and $\gamma$.

It follows from Lemma \ref{selfadjoint} that the operators $\{\I_n(\gamma)\}_{n\in\mathbb{N}}$ are self-adjoint with respect to the inner product (\ref{inner}). By using (\ref{cond2}) and (\ref{Inm}), we then obtain that
\begin{align}
	(\I_n(\gamma)P_{\lambda}, m_{\mu})&=\sum_{\nu\in \mathcal{P}_N, \nu\le_D\lambda}c^{(n)}_{\lambda,\nu}(P_{\nu},m_{\mu})\nonumber\\
	&=(P_{\lambda}, \I_n(\gamma)m_{\mu})=\sum_{\nu\in\mathcal{P}_N,\nu\le_D\mu}b^{(n)}_{\mu,\nu}(P_{\lambda},m_{\nu})=0,
\end{align}
for any $\mu<_D \lambda$.
It implies that $c^{(n)}_{\lambda,\nu}=0$ unless $\lambda=\nu$. Let $G_{n,\lambda}(\gamma)=c^{(n)}_{\lambda,\lambda}$, then we have
\begin{align}\label{Indiag}
	\I_n(\gamma)P_{\lambda}(\x)=G_{n,\lambda}(\gamma)P_{\lambda}(\x).
\end{align}
Thus we finish our proof.
\end{proof}

It follows from (\ref{Indiag}) that
\begin{align}
	[\I_n(\gamma),\I_m(\gamma)]P_{\lambda}(\x)=0
\end{align} 
for any $n,m\in\mathbb{N}$ and $\lambda\in\mathcal{P}_N$. Since $P_{\lambda}(\x)$ is a basis of $\mathbb{C}[\x]^{\mathfrak{S}_N}$, we then obtain
\begin{coro}
The operators $\I_{n}(\gamma)$ in (\ref{Inc}) commute with each other, i.e.,
\begin{align}
	[\I_n(\gamma),\I_m(\gamma)]=0, \quad \forall n,m\in\mathbb{N}.
\end{align}
\end{coro}

\subsection{The eignvalues of $\I_n(\gamma)$}
Taking a basis function $g\in\C[\x]$, we denote $\mathbf{coeff}[g]f$ the coefficient of $g(\x)$ in the expansion of $f(\x)\in\C[\x]$. 
We define a set of $N$-dimensional column vectors
\begin{align}\label{kn1N}
	\kappa^{(n)}(\gamma)=(\kappa_1^{(n)}(\gamma),\cdots,\kappa_N^{(n)}(\gamma))^T
\end{align}
for $n\in\mathbb{N}$, where the elements are given by
\begin{align}\label{kni}
	\kappa_i^{(n)}(\gamma):=\mathbf{coeff}[x^{\lambda}](\pi_{i}^n(\gamma)P_{\lambda}), \quad 1\le i\le N.
\end{align}
Since $\mathbf{coeff}[m_{\lambda}]P_{\lambda}=\mathbf{coeff}[x^{\lambda}]m_{\lambda}=1$, we have $\kappa^{(0)}(\gamma)=\ee=(1,\cdots,1)^T$.

\begin{lemm}
The vectors $\kappa^{(n)}(\gamma)$ in (\ref{kn1N}) satisfy the recurrence formula as
\begin{align}\label{Recur}
	\kappa^{(n+1)}(\gamma)=M_\lambda\kappa^{(n)}(\gamma), \quad \forall n\in\mathbb{N},
\end{align}
where $M_\lambda$ is the $N\times N$ matrix with 
\begin{align}\label{Matrix}
	(M_\lambda)_{ij}=\begin{cases}
		-\beta(\gamma-1) & j<i\\
		\lambda_i-\beta(i-1)+\frac{\beta}{2}(1+\gamma)(N-1) & j=i\\
		-\beta(\gamma+1)& j>i
	\end{cases}.
\end{align}
\end{lemm}

\begin{proof}
It follows from Lemma \ref{selfadjoint} that the operators $\{\pi_{i}(\gamma);i=1,\cdots,N\}$ are covariance, i.e.
\begin{align}
	s_{ij}\pi_{i}(\gamma)s_{ij}=\pi_{j}(\gamma), \quad 1\le i,j\le N.
\end{align}
Since the Jack polynomials are invariant under the action of $\mathfrak{S}_N$, then we have
\begin{align}
	\kappa_j^{(n)}(\gamma)=\mathbf{coeff}[x^{\lambda}](\pi_{j}^n(\gamma)P_{\lambda})=\mathbf{coeff}[x^{s_{ij}\lambda}](\pi_{i}^n(\gamma)P_{\lambda}).
\end{align}
It means that $\kappa_j^{(n)}(\gamma)=\kappa_i^{(n)}(\gamma)$ when $\lambda_i=\lambda_j$.
By (\ref{action}), we have that for any $f\in\C[\x]$ and $\lambda\in \mathcal{P}_N$, the term $x^{\lambda}$ in $\pi_{i}(\gamma)f(\mathbf{x})$ only comes from $\pi_{i}(\gamma)x^{\lambda}$ or $ \pi_{i}(\gamma)x^{s_{ik}\lambda}$ with $k\neq i$ and
\begin{subequations}\label{coeff}
\begin{align}
	\mathbf{coeff}[x^{\lambda}](\pi_{i}(\gamma)x^{\lambda})&=\mathbf{coeff}[x^{\lambda}]\left(x_i\partial_ix^{\lambda}+\frac{\beta}{2}\sum_{j\ne i}\left(\frac{x_i+x_j}{x_i-x_j}+\gamma\right)(x^{\lambda}-x^{s_{ij}\lambda})\right)\nonumber\\
	&=\lambda_i+\frac{\beta}{2}\left(\sum_{j=1,\lambda_i>\lambda_j}^N(\gamma+1)+\sum_{j=1,\lambda_i<\lambda_j}^N(\gamma-1)\right),\\
	\mathbf{coeff}[x^{\lambda}](\pi_{i}(\gamma)x^{s_{ik}\lambda})&=\mathbf{coeff}[x^{\lambda}]\left(x_i\partial_ix^{s_{ik}\lambda}+\frac{\beta}{2}\sum_{j\ne i}\left(\frac{x_i+x_j}{x_i-x_j}+\gamma\right)(x^{s_{ik}\lambda}-x^{s_{ij}s_{ik}\lambda})\right)\nonumber\\
	&=\lambda_i\delta_{\lambda_i,\lambda_k}+\frac{\beta}{2}\left(\sum_{j=1,\lambda_i>\lambda_j}^N(\gamma+1)+\sum_{j=1,\lambda_i<\lambda_j}^N(\gamma-1)\right)\delta_{\lambda_i,\lambda_k}\nonumber\\
	&\quad -\frac{\beta}{2}\left(\sum_{j=1,\lambda_i>\lambda_j}^N(\gamma+1)\delta_{j,k}+\sum_{j=1,\lambda_i<\lambda_j}^N(\gamma-1)\delta_{j,k}\right).
\end{align}
\end{subequations}
By comparing the coefficient of $x^{\lambda}$ on the both sides of the identity $\pi_{i}^{n+1}(\gamma)P_{\lambda}
=\pi_{i}(\gamma)(\pi_{i}^n(\gamma)P_{\lambda})$ and using (\ref{coeff}), we then obtain the recurrence formula
\begin{align}\label{recurrence}
	\kappa_i^{(n+1)}(\gamma)&=\mathbf{coeff}[x^{\lambda}]\left(\kappa_i^{(n)}(\gamma)\pi_{i}(\gamma)x^{\lambda}+\sum_{k=1,\lambda_k\neq \lambda_i}^N\kappa_k^{(n)}(\gamma)\pi_{i}(\gamma)x^{s_{ik}\lambda}\right)\nonumber\\
	&=\lambda_i\kappa_i^{(n)}(\gamma)+\frac{\beta}{2}(\gamma+1)\sum_{j=1,\lambda_i>\lambda_j}^N\kappa_i^{(n)}(\gamma)
	 +\frac{\beta}{2}(\gamma-1)\sum_{j=1,\lambda_i<\lambda_j}^N\kappa_i^{(n)}(\gamma)\nonumber\\
	&\quad -\frac{\beta}{2}(\gamma+1)\sum_{j=1,\lambda_i>\lambda_j}^N\kappa_j^{(n)}(\gamma) -\frac{\beta}{2}(\gamma-1)\sum_{j=1, \lambda_i<\lambda_j}^N\kappa_j^{(n)}(\gamma)\nonumber\\
	&=\lambda_i\kappa_i^{(n)}(\gamma)+\frac\beta2(\gamma+1)\sum_{j>i}(\kappa_i^{(n)}(\gamma)-\kappa_j^{(n)}(\gamma))
	+\frac\beta2(\gamma-1)\sum_{j<i}(\kappa_i^{(n)}(\gamma)-\kappa_j^{(n)}(\gamma))\nonumber\\
	&=\rho_i(\gamma)\kappa_i^{(n)}(\gamma)-\frac\beta2(\gamma+1)\sum_{j>i}\kappa_j^{(n)}(\gamma)
	-\frac\beta2(\gamma-1)\sum_{j<i}\kappa_j^{(n)}(\gamma),
\end{align}
where
\begin{align}\label{rhoi}
	\rho_i(\gamma)=\lambda_i-\beta(i-1)+\frac{\beta}{2}(1+\gamma)(N-1).
\end{align}
This completes the proof since (\ref{Recur}) and (\ref{recurrence}) are equivalent.
\end{proof}

\begin{theo}\label{Eigenvalue}
For the series
\begin{align}\label{Gz}
	G_{\lambda}(z;\gamma)&:=\sum_{n=0}^{\infty}G_{n,\lambda}(\gamma)z^{-1-n},
\end{align}
we have
\begin{align}
	G_{\lambda}(z;\gamma)&=\frac{2}{\beta}\frac{1-\prod_{i=1}^N\frac{z-\rho_i(\gamma)-\frac{\beta}{2}(\gamma+1)}{z-\rho_i(\gamma)-\frac{\beta}{2}(\gamma-1)}}{(\gamma+1)-(\gamma-1)\prod_{i=1}^N\frac{z-\rho_i(\gamma)-\frac{\beta}{2}(\gamma+1)}{z-\rho_i(\gamma)-\frac{\beta}{2}(\gamma-1)}}\nonumber\\
	&=\frac{2}{\beta}\frac{1-\prod_{i=1}^N\frac{z-\rho_i(\gamma)-\frac{\beta}{2}(\gamma-1)}{z-\rho_i(\gamma)-\frac{\beta}{2}(\gamma+1)}}{(\gamma-1)-(\gamma+1)\prod_{i=1}^N\frac{z-\rho_i(\gamma)-\frac{\beta}{2}(\gamma-1)}{z-\rho_i(\gamma)-\frac{\beta}{2}(\gamma+1)}},
\end{align} 
where $\rho_i(\gamma)$ is given by (\ref{rhoi}).
\end{theo}
\begin{proof}
It follows from (\ref{Indiag}) and (\ref{kni}) that
\begin{align}\label{eigen2}
	G_{n,\lambda}(\gamma)=\mathbf{coeff}[x^{\lambda}](\I_n(\gamma)P_{\lambda})=\sum_{i=1}^N\kappa_i^{(n)}(\gamma)=\mathbf{e}^T M_\lambda^n \mathbf{e}.
\end{align}
Thus, we can rewrite (\ref{Gz}) by
\begin{align}\label{Gz2}
	G_{\lambda}(z;\gamma)=\mathbf{e}^T (zI_N-M_\lambda)^{-1} \mathbf{e}.
\end{align}
For convenience, we denote the three $N\times N$ matrices
\begin{align}
	A&:=zI_N-M_{\lambda},\nonumber\\
	A_-&:=A-\frac{\beta}{2}(\gamma-1)\mathbf{e}\mathbf{e}^T,\nonumber\\
	A_+&:=A-\frac{\beta}{2}(\gamma+1)\mathbf{e}\mathbf{e}^T.
\end{align}
It follows from (\ref{Matrix}) that $A_-$ and $A_+$ are lower-triangular and upper-triangular matrices, respectively, and their diagonal elements are
\begin{align}
	(A_-)_{ii}&=z-\rho_i(\gamma)-\frac{\beta}{2}(\gamma-1),\nonumber\\
	(A_+)_{ii}&=z-\rho_i(\gamma)-\frac{\beta}{2}(\gamma+1),
\end{align}
for $i=1,2,\cdots,N$. It gives that
\begin{align}\label{detA-+1}
	\det(A_-)=\prod_{i=1}^N\left[z-\rho_i(\gamma)-\frac{\beta}{2}(\gamma-1)\right],\nonumber\\
	\det(A_+)=\prod_{i=1}^N\left[z-\rho_i(\gamma)-\frac{\beta}{2}(\gamma+1)\right].
\end{align}
Then we have
\begin{align}\label{detA-+2}
	\det(A-\frac{\beta}{2}(\gamma\pm 1)\mathbf{e}\mathbf{e}^T)
	=&\det(A)\det(I_N-\frac{\beta}{2}(\gamma \pm 1)A^{-1}\mathbf{e}\mathbf{e}^T)\nonumber\\
	=&\det(A)[1-\frac{\beta}{2}(\gamma \pm 1)\mathbf{e}^TA^{-1}\mathbf{e}]\nonumber\\
	=&\det(A)[1-\frac{\beta}{2}(\gamma \pm 1)G_{\lambda}(z;\gamma)].
\end{align}
Here we used $\det(I_N+uv^T)=1+\sum_{i=1}^N u_iv_i=1+v^Tu$ for $u=A^{-1}\ee$ and $v=\ee$.
By solving the linear system of two quadratic equations (\ref{detA-+2}) for $\det(A)$ and $G_{\lambda}(z;\gamma)$, we obtain 
\begin{align}\label{GA+A-}
	G_{\lambda}(z;\gamma)&=\frac{2}{\beta}\frac{\det(A_-)-\det(A_+)}{(\gamma+1)\det(A_-)-(\gamma-1)\det(A_+)}\nonumber\\
	&=\frac{2}{\beta}\frac{1-\det(A_+)/\det(A_-)}{(\gamma+1)-(\gamma-1)\det(A_+)/\det(A_-)},\nonumber\\
	&=\frac{2}{\beta}\frac{1-\det(A_-)/\det(A_+)}{(\gamma-1)-(\gamma+1)\det(A_-)/\det(A_+)},
\end{align}
where 
\begin{align}
	\frac{\det(A_+)}{\det(A_-)}&=\prod_{i=1}^N\left(1-\frac{\beta}{z-\rho_i(\gamma)-\frac{\beta}{2}(\gamma-1)}\right)\nonumber\\
	&=\prod_{i=1}^N\left(1+\frac{\beta}{z-\rho_i(\gamma)-\frac{\beta}{2}(\gamma+1)}\right)^{-1}.
\end{align}
This completes the proof of Theorem \ref{Eigenvalue}.
\end{proof}

\begin{ex}
Let $c_i=c_i(\lambda)=\lambda_i-(i-1)\beta$, we have
\begin{subequations}\label{G-1G+1}
\begin{align}
	G_{\lambda}(z;1)&=\beta^{-1}\left[1-\prod_{i=1}^N\left(1-\frac{\beta}{z-c_i+(N-1)\beta}\right)\right]\nonumber\\
	&=\sum_{n=0}^{\infty}\frac{1}{z-c_i+(N-1)\beta}\prod_{\substack{j=1 \\ j\ne i}}^{N}\frac{c_j-c_i+\beta}{c_j-c_i},\\
	G_{\lambda}(z;-1)&=\beta^{-1}\left[\prod_{i=1}^N\left(1+\frac{\beta}{z-c_i}\right)-1\right]\nonumber\\
	&=\sum_{n=0}^{\infty}\frac{1}{z-c_i}\prod_{\substack{j=1 \\ j\ne i}}^{N}\frac{c_j-c_i-\beta}{c_j-c_i}.
\end{align}
\end{subequations}
\end{ex}

\begin{rem}
Sergeev and Veselov \cite{Sergeev15} proved commutativity of the operator set $\{\sum_{i=1}^N\pi_i^n\}_{n\in\mathbb{N}}$ and provided the relationship between $\I_n(0)$ and $\I_n(1)$. Dunkl and Gorin \cite{Dunkl25} obtained the eigenvalues $G_{\lambda}(z;\gamma)$ for $\gamma=1$ and introduced non-symmetric joint eigenfunctions for the commuting operator set $\{\sum_{i=1}^N\pi_i^n\}_{n\in\mathbb{N}}$. However, the relationship between these eigenfunctions and non-symmetric Jack polynomials is unknown.
\end{rem}

\section{Non-symmetric $W$-operators}
In this section, we intend to construct families of non-symmetric $W$-operators. They have a delicate algebraic structure and can be realized by the Dunkl operators.
\subsection{Non-symmetric $W$-operators}
\begin{defin}
Let $a\in\C$. By using the Cherednik-Dunkl operators in (\ref{CHPdunkl}), we define
\begin{align}\label{nonW0a}
	\W_0(a)&:=\frac{1}{2}\sum_{i=1}^Nd_i(d_i+2a-2(i-1)\beta-1)\nonumber\\
	&=\frac12\sum_{i=1}^{N}x_i^2\partial_i^2
	+a\sum_{i=1}^{N}x_i\partial_i+\beta\sum_{1\le i< j\le N}
	\frac{x_ix_j}{x_i-x_j}(\partial_i-\partial_j)\nonumber\\
	&\quad -\beta\sum_{1\le i< j\le N}\frac{x_ix_j}{(x_i-x_j)^2}(1-s_{ij}).
\end{align}
\end{defin}
Form (\ref{dinonJ}), we have that the operator $\W_0(a;\x)$ acts diagonally on the non-symmetric Jack polynomials as
\begin{align}
	\W_0(a)E_{\eta}(\x)=\sum_{(i,j)\in\eta^+}[a+(j-1)+(1-i)\beta]E_{\eta}(\x).
\end{align}
By a straightforward calculation, we arrive at
\begin{prop}\label{W0com}
Let $\tilde{N}=N-1+\beta^{-1}$. The operator $\W_0(\tilde{N}\beta;x)$ satisfies
\begin{subequations}
\begin{align}
	[x_i^{-1}, \W_0(\tilde{N}\beta)]&=x_i^{-1}\pi_i,\\
	[\pi_i, \W_0(\tilde{N}\beta)]&=0,
\end{align}
\end{subequations}
where $\pi_i$ is the Polychronakos' Dunkl operator in (\ref{CHPdunkl}).
\end{prop}

\begin{prop}\label{xrcom}
For $r\in\mathbb{N}$ and $i,j=1,\cdots,N$,  we have 
\begin{subequations}
\begin{align}
	&[\pi_i,x_j^{-1}\pi_j^{r}]=\begin{cases}
		-x_j^{-1}\pi_{j}^{r}-\beta\sum_{k\neq j}x_k^{-1}\pi_{k}^{r}s_{jk} & j=i,\\
		\beta x_j^{-1}\pi_{j}^{r}s_{ij}& j\neq i,
	\end{cases}\label{1xr}\\
	&[x_i^{-1}\pi_i^{r},x_j^{-1}\pi_j^{r}]=0. \label{xrxr}
\end{align}	
\end{subequations}
\end{prop}
\begin{proof}
From (\ref{CHPdunkl}), it is not hard to check that
\begin{align}\label{picom}
	[\pi_i,\pi_j]=\beta(\pi_j-\pi_i)s_{ij},\quad   s_{ij}\pi_i=\pi_js_{ij}, \quad i,j=1,\cdots,N,
\end{align}
and the equation (\ref{1xr}) holds for $r=0$. Assuming (\ref{1xr}) holds for a given $r$, then we have
\begin{align}
	[\pi_i,x_j^{-1}\pi_j^{r+1}]&=[\pi_i,x_j^{-1}\pi_j^{r}]\pi_j+x_j^{-1}\pi_j^{r}[\pi_i,\pi_j]\nonumber\\
	&=\begin{cases}
		\left(-x_j^{-1}\pi_{j}^{r}-\beta\sum_{k\neq j}x_k^{-1}\pi_{k}^{r}s_{jk}\right)\pi_j & j=i,\\
		\beta x_j^{-1}\pi_{j}^{r}[s_{ij}\pi_{j}+(\pi_j-\pi_i)s_{ij}]& j\neq i
	\end{cases}\nonumber\\
	&=\begin{cases}
		-x_j^{-1}\pi_{j}^{r+1}-\beta\sum_{k\neq j}x_k^{-1}\pi_{k}^{r+1}s_{jk} & j=i,\\
		\beta x_j^{-1}\pi_{j}^{r+1}s_{ij}& j\neq i.
	\end{cases}
\end{align}
Thus, we have (\ref{1xr}) holds for any $r\in \mathbb{N}$ by induction. 

Let us turn to proving (\ref{xrxr}).
Since (\ref{xrxr}) holds for the case of $r=0$ or $i=j$, we only need to prove (\ref{xrxr}) for $i\neq j$. Assuming (\ref{xrxr}) holds for a given $r$, we have
\begin{align}
&[x_i^{-1}\pi_i^{r+1},x_j^{-1}\pi_j^{r+1}]\nonumber\\
=&x_i^{-1}\pi_i^{r}[\pi_i,x_j^{-1}\pi_j^{r}]\pi_j-x_j^{-1}\pi_j^{r}[\pi_j,x_i^{-1}\pi_i^{r}]\pi_i +x_i^{-1}\pi_i^{r}x_j^{-1}\pi_j^{r}[\pi_i,\pi_j]\nonumber\\
=&\beta x_i^{-1}\pi_i^{r}x_j^{-1}\pi_{j}^{r}s_{ij}\pi_j-\beta x_j^{-1}\pi_j^{r}x_i^{-1}\pi_{i}^{r}s_{ij}\pi_i +\beta x_i^{-1}\pi_i^{r}x_j^{-1}\pi_j^{r}(\pi_j-\pi_i)s_{ij}\nonumber\\
=&\beta x_i^{-1}\pi_i^{r}x_j^{-1}\pi_{j}^{r}[s_{ij}\pi_{j}-s_{ij}\pi_i+(\pi_j-\pi_i)s_{ij}]=0.
\end{align}
Thus, we proved (\ref{xrxr}) by induction.
\end{proof}

The Proposition \ref{W0com} and Proposition \ref{xrcom} suggest that the operators $\{x_i^{-1}\pi_i^r;r\in\mathbb{N},i=1,\cdots,N\}$ and $\W_0(a)$ generate an algebra that is expected to be a subalgebra of dDAHA of $GL(N)$. In order to investigate the algebra more clearly, we will introduce a set of ``non-symmetric" $W$-operators associated with $\pi(\gamma)$ in (\ref{HPDunkl}).
\begin{defin}\label{nonWdef}
Let $r\in\mathbb{N}$, $k\in\mathbb{Z}_{>0}$, $a_0\in\C$ and $\mathbf{a}=(a_1,\cdots,a_r)\in\C^r$, we define a set of $W$-operators on the space $\C[\x^{\pm}]$
\begin{align}\label{nonWop}
	\W_{\pm k,i}^{(r)}(\mathbf{a})=\begin{cases}
		\prod_{l=1}^r\mathbf{ad}_{\pm \W_0(a_l)}x_i^{\pm 1} & k=1,\\
		\frac{1}{(k-1)!}\mathbf{ad}^{k-1}_{\pm \sum_{j=1}^N\W_{\pm 1,j}^{(r+1)}(a_0,\mathbf{a})}\W_{\pm 1,i}^{(r)}(\mathbf{a}) & k\ge 1,
	\end{cases}
\end{align}
for $i=1,\cdots,N$, where $\mathbf{ad}_X Y=[X,Y]$. These $W$-operators are independent to the value of $a_0$.
\end{defin}

We call $\{\W_{-k,i}^{(r)}(\mathbf{a})\}$ and $\{\W_{k,i}^{(r)}(\mathbf{a})\}$ the negative and positive branches, respectively. It follows from Lemma \ref{selfadjoint} that $\W_0(a)$ is self-adjoint with respect to the inner product (\ref{noninner}).
Then we have
\begin{align}\label{nonW-11}
	(\W_{-1,i}^{(r+1)}(\tilde{N}\beta,\mathbf{a}))^{\dagger}&=(\prod_{l=1}^r\mathbf{ad}_{-\W_0(a_l)}\D_i)^\dagger\nonumber\\
	&=\prod_{l=1}^r\mathbf{ad}_{\W_0(a_l)}x_i=\W_{1,i}^{(r)}(\mathbf{a}).
\end{align}
Taking the action of ``$\dagger$" on the both sides of (\ref{nonWop}), we further obtain
\begin{prop}\label{W+W-}
There is a duality between the $W$-operators in negative and positive branches
\begin{align}
	(\W_{-k,i}^{(r+1)}(\tilde{N}\beta,\mathbf{a}))^{\dagger}=\W_{k,i}^{(r)}(\mathbf{a}).
\end{align}
\end{prop}

\subsection{Factorization by Dunkl operators}
In order to give more direct representation for $\{\W_{\pm k,i}^{(r)}(\mathbf{a};\x)\}$, we define the following $W$-operators degree $0$
\begin{subequations}\label{W0r}
\begin{align}
	\W_{0,i}^{(r)}(\mathbf{a})&=\prod_{l=1}^r(\pi_i+a_l-\tilde{N}\beta),\\
	\tilde{\W}_{0,i}^{(r)}(\mathbf{a})&=\prod_{l=1}^r(\tilde{\pi}_i+a_l),
\end{align}
\end{subequations}
where $\tilde{\pi}_i=\tilde{\pi}_i(\x)=\pi(-1;\x)$. For $\tilde{\pi}_i$, we have
\begin{align}\label{tildepi}
	\tilde{\pi}_i(\x)+\tilde{N}\beta&=x_i^{-1}\circ\pi_i(\x)\circ x_i= \D_ix_i,\nonumber\\
	\tilde{\pi}_i(\x)&=-\pi_i(\x^{-1}).
\end{align}

Now, we propose the main result in this section.
\begin{theo}\label{Woprep}
The $W$-operators in Definition \ref{nonWdef} can be factored by the Dunkl operators as
\begin{subequations}
\begin{align}
	\W_{-k,i}^{(r)}(\mathbf{a})&=\left[x_i^{-1}\W_{0,i}^{(r)}(\mathbf{a})\right]^k,\label{W-k}\\
	\W_{k,i}^{(r)}(\mathbf{a})&=\left[x_i\tilde{\W}_{0,i}^{(r)}(\mathbf{a})\right]^k.\label{W+k}
\end{align}
\end{subequations}
\end{theo}

\begin{proof}
We denote Euler operator $E=\sum_{i=1}^Nx_i\partial_i$.
By Proposition \ref{W0com}, we then get that
\begin{align}
	\W_{-1,i}^{(r)}(\mathbf{a})&=\prod_{l=1}^r[\mathbf{ad}_{-\W_0(\tilde{N\beta})}-(a_l-\tilde{N}\beta)\mathbf{ad}_{E}] x_i^{-1} \nonumber\\
	&=x_i^{-1}\prod_{l=1}^r(\pi_i+a_l-\tilde{N}\beta)=x_i^{-1}\W_{0,i}^{(r)}(\mathbf{a}).
\end{align}
It is known that $\W_{0,i}^{(r)}(\mathbf{a})=\sum_{l=1}^re_{r-n}(\mathbf{a})\pi_i^{n}$.
By using (\ref{1xr}) and (\ref{xrxr}), we then obtain
\begin{align}
    &[\pi_i,\W_{-1,j}^{(r)}(\mathbf{a})]=\begin{cases}
	-\W_{-1,j}^{(r)}(\mathbf{a})-\beta\sum_{k\neq j}\W_{-1,k}^{(r)}(\mathbf{a})s_{jk} & j=i,\\
	\beta \W_{-1,j}^{(r)}(\mathbf{a})s_{ij}& j\neq i,
    \end{cases}\\
	&[\W_{-1,i}^{(r)}(\mathbf{a}),\W_{-1,j}^{(r)}(\mathbf{a})]=0. \label{Wij}
\end{align}
It further gives 
\begin{align}
	&\sum_{j=1}^N[\W_{-1,j}^{(r+1)}(a_0,\mathbf{a}),\W_{-1,i}^{(r)}(\mathbf{a})]\nonumber\\
	=&\sum_{j=1}^N\W_{-1,j}^{(r)}(\mathbf{a})[a_0+\pi_j,\W_{-1,i}^{(r)}(\mathbf{a})]+\sum_{j=1}^N[\W_{-1,j}^{(r)}(\mathbf{a}),\W_{-1,i}^{(r)}(\mathbf{a})](a_0+\pi_j)\nonumber\\
	=&-(\W_{-1,i}^{(r)}(\mathbf{a}))^2.
\end{align}
Therefore, we obtain
\begin{align}
	\W_{-k,i}^{(r)}(\mathbf{a})&=\frac{1}{(k-1)!}\mathbf{ad}^{k-1}_{-\sum_{j=1}^N\W_{-1,j}^{(r+1)}(a_0,\mathbf{a})} \W_{-1,i}^{(r)}(\mathbf{a;x})\nonumber\\
	&=\frac{1}{(k-1)!}\mathbf{ad}^{k-2}_{-\sum_{j=1}^N\W_{-1,j}^{(r+1)}(a_0,\mathbf{a})} (\W_{-1,i}^{(r)}(\mathbf{a}))^2\nonumber\\
	&=\left(\W_{-1,i}^{(r)}(\mathbf{a})\right)^{k}.
\end{align}
Here we used $[X,Y^n]=\sum_{k=0}^{n-1}Y^{k}[X,Y]Y^{n-k-1}$.

It follows form Proposition \ref{W+W-} that
\begin{align}
	\W_{k,i}^{(r)}(\mathbf{a})&=\left[\left(x_i^{-1}\pi_i\prod_{l=1}^r(\pi_i+a_l-\tilde{N}\beta)\right)^k\right]^\dagger\nonumber\\
	&=\left[\left(\prod_{l=1}^r(\tilde{\pi}_i+a_l)\D_i\right)^\dagger\right]^k=\left[x_i\prod_{l=1}^r(\tilde{\pi}_i+a_l)\right]^k,
\end{align}
where we used (\ref{tildepi}) and self-adjointness of $\pi_i$ and $\tilde{\pi}_i$.
Thus we finish our proof.
\end{proof}

Follows from Equations (\ref{tildepi}), (\ref{Wij}) and Theorem \ref{Woprep}, we conclude that 
\begin{coro}\label{Wdual}
There is another duality between the operators in the negative and positive branches
\begin{align}
	\W_{-k,i}^{(r)}(\tilde{N}\beta-\mathbf{a};\x)=(-1)^{kr}\W_{k,i}^{(r)}(\mathbf{a};\x^{-1}),
\end{align}
where $\tilde{N}\beta-\mathbf{a}=(\tilde{N}\beta-a_1,\cdots,\tilde{N}\beta-a_N)$.
\end{coro}
\begin{coro}\label{noncom}
The $W$-operators in Definition \ref{nonWdef} generate the following commuting algebra
\begin{subequations}
\begin{align}
	[\W_{\pm k,i}^{(r)}(\mathbf{a}),\W_{\pm k,j}^{(r)}(\mathbf{a})]&=0,\\
	[\W_{\pm k,i}^{(r)}(\mathbf{a}),\W_{\pm k,i}^{(r')}(\mathbf{a})]&=0,
\end{align}
\end{subequations}
for any $i,j=1,\cdots,N$, $r,r'\in\mathbb{N}$.
\end{coro}

\subsection{Representations by auxiliary operator}
Let us introduce a auxiliary operator $\OO(a)$ which satisfies
\begin{align}\label{nonOop}
	\OO(a)E_{\eta}(\x)=[a]_{\eta_+}E_{\eta}(\x),\quad \forall \eta\in V_N.
\end{align}
It is clear that $\OO(a)$ is invertible on $\C[\x]$ when $[a]_{\eta^+}\neq 0$. By (\ref{dinonJ}) and (\ref{SnonJ}), it gives that 
the operators $\{d_i\}_{i=1,\cdots,N}$ and $w\in \mathfrak{S}_N$ respectively commute with  $\OO(a;\x)$. It further gives 
\begin{align}
	[\OO(a),\W_{0,i}^{(r)}(\mathbf{a})]=[\OO(a),\tilde{\W}_{0,i}^{(r)}(\mathbf{a})]=0, \ i=1,\cdots,N.
\end{align}

\begin{theo}\label{nonWO}
Let $\mathbf{a;x}=(a_1,\cdots,a_r)$ with $a_l\in\C\setminus\{m+n\beta|(m,n)\in\mathbb{Z}_{\le0}^2\}$ for $l=1,\cdots,r$ and $\OO^{(r)}(\mathbf{a;x})=\prod_{l=1}^r\OO(a_l)$. Then we have
\begin{subequations}
	\begin{align}
		\W_{k,i}^{(r)}(\mathbf{a})&=\mathbf{Ad}_{\OO^{(r)}(\mathbf{a})} \W_{k,i}^{(0)},\\
		\W_{-k,i}^{(r)}(\mathbf{a})&=\mathbf{Ad}^{-1}_{\OO^{(r)}(\mathbf{a})} \W_{-k,i}^{(0)},
	\end{align}
\end{subequations}
for $i=1,\cdots,N$ and $k\in\mathbb{Z}_{>0}$, where the operation $\mathbf{Ad}_X Y=X\circ Y\circ X^{-1}$.
\end{theo}

\begin{proof}
The Theorem 3.3 in Ref.\cite{Brennecken22} shows that
\begin{align}\label{ED}
	E_{\eta}(\D_1(\y),\cdots,\D_N(\y))\prod_{i=1}^Ny_i^{-a}=(-1)^{|\eta|}[a]_{\eta_+}E_{\eta}(\y^{-1})\prod_{i=1}^Ny_i^{-a},
\end{align}
where $\{\D_i\}_{i=1,\cdots,N}$ are the Dunkl operators (\ref{Dunkl}). It is known that $\{E_{\eta}(\x)\}_{\eta\in V_N}$ is a basis of $\C[\x]$. Let $x_i=y_i^{-1}$ for $i=1,\cdots,N$. For any homogeneous polynomial $f(\x)\in\C[\x]$ degree $m$, we have
\begin{align}\label{fD}
	f(\D_1(\y),\cdots,\D_N(\y))\left(\prod_{i=1}^Nx_i^{a}\right)=(-1)^{m}\prod_{i=1}^Nx_i^{a}[\OO(a)f(\x)].
\end{align} 
In addition, a straightforward calculation shows that
\begin{align}
	\D_i(\y)\circ\prod_{i=1}^Nx_i^{a}=-\prod_{i=1}^Nx_i^{a}\circ[x_i(\tilde{\pi}_i+a)].
\end{align} 
Taking the action of $\D_i(\y)$ on the both sides of (\ref{fD}), it gives that
\begin{align}
    f(\D_1(\y),\cdots,\D_N(\y))\left(\D_i(\y)\prod_{i=1}^Nx_i^{a}\right)
    &=(-1)^{m}\D_i(\y)\left(\prod_{i=1}^Nx_i^{a}\OO(a)f(\x)\right)\nonumber\\
    &=(-1)^{m+1}\prod_{i=1}^Nx_i^{a}[x_i(\tilde{\pi}_i+a)]\OO(a)f(\x)\nonumber\\
    &=(-1)^{m+1}\prod_{i=1}^Nx_i^{a}[\OO(a)(x_if)(\x)],
\end{align}
where we used (\ref{fD}) by taking replace $f(\x)$ with $x_if(\x)$ in the last line. Thus we obtain that
\begin{align}\label{OxxO}
    \OO(a)\circ x_i=[x_i(\tilde{\pi}_i+a)]\circ \OO(a).
\end{align}
Since the non-symmetric Jack polynomial $E_{\eta}$ is orthogonal with respect to the inner product (\ref{noninner}) and $\OO(a)$ acts diagonally on it, then $\OO(a)$ is self-adjoint. By taking the adjoint map ``$\dagger$" (\ref{adjoint2}) on the both sides of (\ref{OxxO}), we get that
\begin{align}\label{xOOx}
	x_i^{-1}\pi_i\circ\OO(a)&=\OO(a;\x)\circ[(\tilde{\pi}_i+a)\D_i]\nonumber\\
	&=\OO(a)\circ[x_i^{-1}\pi_i(\pi_i+a-\tilde{N}\beta)],
\end{align}
where we used (\ref{tildepi}).
By (\ref{OxxO}) and (\ref{xOOx}), we further obtain that
\begin{subequations}
\begin{align}
	[x_i\prod_{l=1}^r(\tilde{\pi}_i+a_l)]^k&=\mathbf{Ad}_{\OO^{(r)}(\mathbf{a})}x_i^k,\\
	[x_i^{-1}\prod_{l=1}^r(\pi_i+a_l-\tilde{N}\beta)]^k&=\mathbf{Ad}^{-1}_{\OO^{(r)}(\mathbf{a})}x_i^{-k}.
\end{align}
\end{subequations}
Thus, we finish our proof.
\end{proof}

\subsection{Symmetrization of $W$-operators}
\begin{defin}
Based on the non-symmetric operators in (\ref{nonW0a}), (\ref{nonWop}) and (\ref{W0r}), we define families of $W$-operators by following symmetrization
\begin{subequations}\label{Wsymm}
\begin{align}
	W_{0}(\mathbf{a})&:=\Res\W_{0}(\mathbf{a}),\\
	W_{0}^{(r)}(\mathbf{a})&:=\Res\left(\sum_{i=1}^N \W_{0,i}^{(r)}(\mathbf{a})\right),\\
	\tilde{W}_{0}^{(r)}(\mathbf{a})&:=\Res\left(\sum_{i=1}^N \tilde{\W}_{0,i}^{(r)}(\mathbf{a})\right),\\
	W_{\pm k}^{(r)}(\mathbf{a})&:=\Res\left(\sum_{i=1}^N \W_{\pm k,i}^{(r)}(\mathbf{a})\right), \label{Wopdef}
\end{align}
\end{subequations}
for $r\in\mathbb{N}$ and $k\in\mathbb{Z}_{>0}$.
\end{defin}
It follows from Lemma \ref{selfadjoint} that these symmetrized $W$-operators are invariant under the actions of the permutation group $\mathfrak{S}_N$.
At first, we consider these $W$-operators degree $0$. By (\ref{Inc}) and (\ref{W0r}), it gives 
\begin{subequations}
\begin{align}
	W^{(r)}_{0}(\mathbf{a})&=\sum_{n=0}^{r}e_{r-n}(\mathbf{a}-\tilde{N}\beta)\I_n(1),\\
	\tilde{W}^{(r)}_{0}(\mathbf{a})&=\sum_{n=0}^{r}e_{r-n}(\mathbf{a})\I_n(-1).
\end{align}
\end{subequations}
Follows from Theorem \ref{Eigenvalue} and Equation (\ref{G-1G+1}), we then get
\begin{prop}\label{W0rs}
The operators $\{W^{(r)}_{0}(\mathbf{a}),\tilde{W}^{(r)}_{0}(\mathbf{a})\}_{r\in\mathbb{N}}$ act diagonally on the Jack polynomials, i.e. 
\begin{subequations}
\begin{align}
	W^{(s)}_{0}(\mathbf{a;x})C_{\lambda}(\x)&=G_\lambda^{(s)}(\mathbf{a})C_{\lambda}(\x),\\
	\tilde{W}^{(r)}_{0}(\mathbf{a;x})C_{\lambda}(\x)&=\tilde{G}_\lambda^{(r)}(\mathbf{a})C_{\lambda}(\x),
\end{align}
\end{subequations}
and the eigenvalues are given by
\begin{subequations}\label{GaGb}
\begin{align}
	G_\lambda^{(r)}(\mathbf{a})&=\sum_{n=0}^re_{r-n}(\mathbf{a}-\tilde{N}\beta)\mathbf{coeff}[z^{-1-n}]G_{\lambda}(z;1)\nonumber\\
	&=\sum_{i=1}^{N}\prod_{l=1}^r(a_l+c_i-1)\prod_{\substack{j=1 \\ j\ne i}}^{N}\frac{c_j-c_i+\beta}{c_j-c_i},\\
	\tilde{G}_\lambda^{(r)}(\mathbf{a})&=\sum_{n=0}^re_{r-n}(\mathbf{a})\mathbf{coeff}[z^{-1-n}]G_{\lambda}(z;-1)\nonumber\\
	&=\sum_{i=1}^{N} \prod_{l=1}^r(a_l+c_i)\prod_{\substack{j=1 \\ j\ne i}}^{N}\frac{c_j-c_i-\beta}{c_j-c_i}.
\end{align}
\end{subequations}
\end{prop}

\begin{ex}
We list the complete form for some $W$-operators in (\ref{Wopdef}) 
\begin{subequations}
\begin{align}
	W_{\pm k}^{(0)}(\x)&=p_{\pm k}(\x),\\
	W_{-1}^{(1)}(a,\x)&=\sum_{i=1}^N\partial_i+(a-\tilde{N}\beta)p_{-1}(\x),\\
	W_{1}^{(1)}(a,\x)&=\sum_{i=1}^Nx_i^2\partial_i+ap_1(\x),\\
	W_{-1}^{(2)}(a,b;\x)&=\sum_{i=1}^{N} [\partial_ix_i\partial_i
	+(a+b-2\tilde{N}\beta) \partial_i ]+ 2\beta \sum_{1\le i\ne j\le N}
	\frac{1}{x_i-x_j}x_i\partial_i\nonumber\\
	&\quad  +(a-\tilde{N}\beta)(b-\tilde{N}\beta)p_{-1}(\x),\\
	W_{1}^{(2)}(a,b;\x)&=\sum_{i=1}^N[x_i(x_i\partial_i)^2
	+(a+b)x_i^2\partial_i]+2\beta\sum_{1\le i\neq j\le N}
	\frac{x_ix_j}{x_i-x_j}x_i\partial_i+abp_1(\x).
\end{align}
\end{subequations} 
\end{ex}
We can see that the rational-function terms in the operators $\{W_{\pm 1}^{(r)}\}$ become increasingly complicated when the parameters $r$ increase. This indicates that the simplification of these $W$-operators through the Dunkl operator is quite necessary. 

Follows from Theorem \ref{Woprep}, Corollary \ref{noncom}, Corollary \ref{Wdual} and Theorem \ref{nonWO}, we conclude that
\begin{coro}\label{Wrep}
For the symmetric operators $\{W_{\pm k}^{(r)}(\mathbf{a;x}):r\in\mathbb{N},\ k\in\mathbb{Z}_{>0}\}$ in (\ref{Wopdef}), we have \\
(1) They can be identified by
\begin{align}\label{Wop}
	W^{(r)}_{\pm k}(\mathbf{a})&=\begin{cases}
		\prod_{l=1}^r\mathbf{ad}_{\pm W_0(a_l)}p_{\pm 1}(\x) & k=1,\\
		\frac{1}{(k-1)!}\mathbf{ad}^{k-1}_{\pm W^{(r+1)}_{\pm 1}(a_0,\mathbf{a})}W^{(r)}_{\pm 1}(\mathbf{a}) &  k>1.
	\end{cases}
\end{align}
(2) They can be rewritten as 
\begin{subequations}\label{sumW-+k}
	\begin{align}
		W_{-k}^{(r)}(\mathbf{a})&=\Res \sum_{i=1}^N\left[x_i^{-1}\prod_{l=1}^r(\pi_i+a_l-\tilde{N}\beta)\right]^k, \label{sumW-k}\\
		W_{k}^{(r)}(\mathbf{a})&=\Res \sum_{i=1}^N \left[x_i\prod_{l=1}^r(\tilde{\pi}_i+a_l)\right]^k.
	\end{align}
\end{subequations}
Especially, $W_{k}^{(r)}(\mathbf{a})$ and $W_{-k}^{(r+1)}(\tilde{N}\beta,\mathbf{a})$ are differential operators on $\C[\x]^{\mathfrak{S}_N}$.\\
(3) They can be represented by the auxiliary operator
\begin{subequations}\label{WO}
	\begin{align}
		W_{k}^{(r)}(\mathbf{a})&=\mathbf{Ad}_{\OO^{(r)}(\mathbf{a})} p_k(\x),\\
		W_{-k}^{(r)}(\mathbf{a})&=\mathbf{Ad}^{-1}_{\OO^{(r)}(\mathbf{a})} p_k(\x^{-1}),
	\end{align}
\end{subequations}
when $a_l\in\C\setminus\{m+n\beta|(m,n)\in\mathbb{Z}_{\le0}^2\}$ for $l=1,\cdots,r$.\\
(4) There is a duality between the operators in the negative and positive branches
\begin{align}
	W_{-k}^{(r)}(\tilde{N}\beta-\mathbf{a};\x)&=(-1)^{kr}\W_{k}^{(r)}(\mathbf{a};\x^{-1}).
\end{align}
(5) They generate a commuting algebra as
\begin{align}\label{Wcomm}
	[W_{\pm k}^{(r)}(\mathbf{a}),W_{\pm k}^{(r')}(\mathbf{a})]=0,\quad \forall r,r'\in\mathbb{N}.
\end{align}
\end{coro}
It is worth noting that the operation ``$\Res$" is linear but not homomorphic. Generally speaking, 
\begin{align}
	\Res X\circ\Res Y\neq \Res(X\circ Y).
\end{align}
unless $Y$ is the operator on $\C[\x^{\pm }]^{\mathfrak{S}_N}$. Thus, (\ref{Wop}) can not be obtained by summing over index $i$ and taking ``$\Res$" on every operators in both sides of (\ref{nonWop}) directly but calculating recursively.

\begin{rem}
Here we present three different representations (\ref{sumW-+k}), (\ref{Wop}) and (\ref{WO}) for the symmetric operators $\{W_{\pm k}^{(r)}(\mathbf{a}):k\in\mathbb{Z}_{>0},r\in\mathbb{N}\}$.
In \cite{Liu24}, we defined them by (\ref{Wop}) but not (\ref{Wopdef}). Some crucial propositions related to these symmetric operators including the commutativity (\ref{Wcomm}) and hypergeometric equations mentioned later were obtained by (\ref{WO}). In addition, (\ref{sumW-k}) for the case of $\mathbf{a}=(\tilde{N}\beta,\cdots,\tilde{N}\beta)$ was first given by \cite{Mironov23} but without a proof. 
We provide a proof for (\ref{sumW-+k}) and (\ref{Wcomm}) by the symmetrization (\ref{Wopdef}) here.
\end{rem}

\section{Macdonald's $\beta$-deformed hypergeometric series}
\subsection{Hypergeometric series}
Let us briefly review the non-symmetric Macdonald's $\beta$-deformed hypergeometric series defined by Brennecken and R\"osler.
\begin{defin}\label{nonMac}\cite{Brennecken22}
Fix $r,s\in\mathbb{N}$ and let $\mathbf{a}=(a_1,\cdots,a_r)\in\C^r$ and $\mathbf{b}=(b_1,\cdots,b_s)\in \C^s$ be such that $b_l\notin\{m+n\beta|(m,n)\in\mathbb{Z}_{\le0}^2\}$. The non-symmetric Macdonald's $\beta$-deformed hypergeometric series are defined by
\begin{align}
	_r\mathcal{K}_s(\mathbf{a;b;x;y})&:=\sum_{\eta\in V_N}\frac{\prod_{k=1}^r[a_k]_{\eta_+}} {\prod_{l=1}^s[b_l]_{\eta_+}}\frac{L_{\eta}(\mathbf{x})L_{\eta}(\mathbf{y})}{L_{\eta}(1^N)|\eta|!},
\end{align}
as well as the symmetric hypergeometric series
\label{Hyper}
\begin{align}
	_r\mathcal{F}_s(\mathbf{a;b;x;y}):=\sum_{\lambda\in\mathcal{P}_N}\frac{\prod_{k=1}^r[a_k]_{\lambda}}{\prod_{l=1}^s[b_l]_{\lambda}}\frac{C_{\lambda}(\mathbf{x})C_{\lambda}(\mathbf{y})}{C_{\lambda}(1^N)|\lambda|!},
\end{align}
where $[a]_{\lambda}$ and $L_{\eta}(\x)$ are respectively given by (\ref{contants}) and (\ref{nonL}).
\end{defin}
The Macdonald's $\beta$-deformed hypergeometric series with single alphabet $\x$ are defined by
\begin{align}\label{Frsx}
	_rF_s(\mathbf{a;b;x})={_r\mathcal{F}_s}(\mathbf{a;b;x};1^N)=\sum_{\lambda\in\mathcal{P}_N}\frac{\prod_{k=1}^r[a_k]_{\lambda}}{\prod_{l=1}^s[b_l]_{\lambda}}\frac{C_{\lambda}(\mathbf{x})}{|\lambda|!}.
\end{align}

\begin{prop}\cite{Brennecken22}\label{FKFrs}
The non-symmetric and symmetric hypergeometric series are connected by
\begin{align}
	_r\mathcal{F}_s(\mathbf{a;b;x;y})=\frac{1}{N!}\sum_{w\in\mathfrak{S}_N}{_r\mathcal{K}_s}(\mathbf{a;b};\x;w\y).
\end{align}
Especially, it gives $_rF_s(\mathbf{a;b;x})={_r\mathcal{K}_s}(\mathbf{a;b;x};1^N)$.
\end{prop}

By using the auxiliary operator $\OO(a)$ in (\ref{nonOop}), we immediately get that
\begin{prop}\label{KrsK00}
The non-symmetric hypergeometric series in Definition \ref{nonMac} can be rewritten as
\begin{align}
	_r\mathcal{K}_s(\mathbf{a;b;x;y})&=\OO^{(r)}(\mathbf{a};\x)\left(\OO^{(s)}(\mathbf{b};\y)\right)^{-1}{_0\mathcal{K}_0}(\x;\y).
\end{align}
\end{prop}

\subsection{Hypergeometric equations}
In this part, we will show that the symmetric and non-symmetric $W$-operators defined in section 4 play a crucial role in characterizing the hypergeometric series. At first, let us consider the series in two sets of variables.
\begin{lemm}\cite{Baker98}\label{xynonK}
The series ${_0\mathcal{K}_0}(\x;\y)$ satisfies the following constraints
\begin{align}\label{K00}
	[x_i^k-(y_i^{-1}\pi_i(\y))^k]{_0\mathcal{K}_0}(\x;\y)=0, \quad i=1,\cdots,N,\ k\in\mathbb{Z}_{>0}.
\end{align} 
\end{lemm}
By taking the actions of $\OO^{(r)}(\mathbf{a};\x)\left(\OO^{(s)}(\mathbf{b};\y)\right)^{-1}$ on the both sides of (\ref{K00}) and applying Theorem \ref{nonWO} and Proposition \ref{KrsK00}, we have
\begin{theo}
The non-symmetric hypergeometric series ${_r\mathcal{K}_s}(\mathbf{a;b};\x;\y)$ can be characterized by the following equation system
\begin{align}\label{nonWK}
	\left(\W_{k,i}^{(r)}(\mathbf{a;x})-\W_{-k,i}^{(s+1)}(\tilde{N}\beta,\mathbf{b;y})\right)F(\x;\y)=0, \quad \forall i=1,\cdots,N, \ k\in\mathbb{Z}_{\ge 0},
\end{align}
\end{theo}
It is known that the symmetric operators $\{W_{k}^{(r)}(\mathbf{a;x}):k\in\mathbb{Z}_{>0}, \ r\in\mathbb{N}\}$ in (\ref{Wop}) are invariant under the action of $\mathfrak{S}_N$ from Corollary \ref{Wrep}. By summing over the index $i$ and then taking the action of $\frac{1}{N!}\sum_{w\in \mathfrak{S}_N}w^{(\y)}$ on (\ref{nonWK}) and applying Proposition \ref{FKFrs}, we recover that
\begin{coro}\label{Wxycons}\cite{Liu25}
The symmetric hypergeometric series ${_r\mathcal{F}_s}(\mathbf{a;b};\x;\y)$ can be characterized by the following equation system
\begin{align}
	\left(W_{k}^{(r)}(\mathbf{a;x})-W_{-k}^{(s+1)}(\tilde{N}\beta,\mathbf{b;y})\right)F(\x;\y)=0, \ k\in\mathbb{Z}_{\ge 0}.
\end{align}
\end{coro} 
Let us turn to the hypergeometric series in a single set of variables.
We briefly review Chen and Sahi’s results associated with their hypergeometric equations. 
\begin{theo}\cite{Chen25}\label{W0W+-}
The hypergeometric series $_rF_s(\mathbf{a;b;x})$ satisfy 
\begin{subequations}
\begin{align}
	\left(\tilde{W}^{(r)}_{0}(\mathbf{a;x})-W^{(s+1)}_{-1}(\tilde{N}\beta,\mathbf{b;x})\right)F(\mathbf{x})&=0, \\
	\left(W^{(r)}_{1}(\mathbf{a;x})-W^{(s+1)}_{0}(\tilde{N}\beta,\mathbf{b;x})\right)F(\mathbf{x})&=0.
\end{align}
\end{subequations}
\end{theo}
\begin{proof}
By the generalized binomial coefficients (\ref{binom}), the eigenvalues $\tilde{G}^{(r)}_\lambda(\mathbf{a})$ and $G^{(s)}_\lambda(\mathbf{b})$ in (\ref{GaGb}) can be rewritten as
\begin{subequations}
\begin{align}
	\tilde{G}^{(r)}_\lambda(\mathbf{a})&=\frac{1}{|\lambda|+1}\sum_{i=1}^N\prod_{k=1}^r(a_k+c_i)\binom{\lambda^{(i)}}{\lambda}_{\beta}\frac{C_{\lambda^{(i)}}(1^N)}{C_{\lambda}(1^N)},\\
	G^{(s+1)}_\lambda(\tilde{N}\beta,\mathbf{b})&=\sum_{i=1}^N\prod_{l=1}^s(b_l+c_i-1)\binom{\lambda}{\lambda_{(i)}}_{\beta}.
\end{align}
\end{subequations}
By using the Pieri formulas (\ref{Pieri}) and Proposition \ref{W0rs}, we get that
\begin{align}
	&W^{(s+1)}_{-1}(\tilde{N}\beta,\mathbf{b;x}){_rF_s}(\mathbf{a;b;x})\nonumber\\
	=&(\OO^{(s)}(\mathbf{b;x}))^{-1}\left(\sum_{i=1}^N\partial_i\right){_rF_0}(\mathbf{a;x})\nonumber\\
	=&\sum_{i=1}^N\sum_{\mu_{(i)}\in\mathcal{P}_N}\frac{\prod_{k=1}^r[a_k]_{\mu}}{\prod_{l=1}^s[b_l]_{\mu_{(i)}}}\binom{\mu}{\mu_{(i)}}_{\beta}\frac{C_{\mu}(1^N)}{C_{\mu_{(i)}}(1^N)}\frac{C_{\mu_{(i)}}(\mathbf{x})}{|\mu|!}\nonumber\\
	=&\sum_{\lambda\in\mathcal{P}_N}\sum_{i=1}^N\frac{\prod_{k=1}^r[a_k]_{\lambda^{(i)}}}{\prod_{l=1}^s[b_l]_{\lambda}}\binom{\lambda^{(i)}}{\lambda}_{\beta}\frac{C_{\lambda^{(i)}}(1^N)}{C_{\lambda}(1^N)}\frac{C_{\lambda}(\mathbf{x})}{(|\lambda|+1)!}\nonumber\\
	=&\sum_{\lambda\in\mathcal{P}_N}\tilde{G}^{(r)}_\lambda(\mathbf{a})\frac{\prod_{k=1}^r[a_k]_{\lambda}}{\prod_{l=1}^s[b_l]_{\lambda}}\frac{C_{\lambda}(\mathbf{x})}{|\lambda|!}\nonumber\\
	=&\tilde{W}^{(r)}_{0}(\mathbf{a;x}){_rF_s}(\mathbf{a;b;x}),
\end{align}
and 
\begin{align}
	&W^{(r)}_{1}(\mathbf{a;x}){_rF_s}(\mathbf{a;b;x})\nonumber\\
	=&\OO^{(r)}(\mathbf{a;x})p_1(\x){_0F_s}(\mathbf{b;x})\nonumber\\
	=&\sum_{i=1}^N\sum_{\mu_{(i)}\in\mathcal{P}_N}\frac{\prod_{k=1}^r[a_k]_{\mu^{(i)}}}{\prod_{l=1}^s[b_l]_{\mu}}\binom{\mu^{(i)}}{\mu}_{\beta}\frac{C_{\mu^{(i)}}(\mathbf{x})}{|\mu^{(i)}|!}\nonumber\\
	=&\sum_{\lambda\in\mathcal{P}_N}\sum_{i=1}^N\frac{\prod_{k=1}^r[a_k]_{\lambda}}{\prod_{l=1}^s[b_l]_{\lambda_{(i)}}}\binom{\lambda}{\lambda_{(i)}}_{\beta}\frac{C_{\lambda}(\mathbf{x})}{|\lambda|!}\nonumber\\
	=&\sum_{\lambda\in\mathcal{P}_N}G^{(s+1)}_\lambda(\tilde{N}\beta,\mathbf{b})\frac{\prod_{k=1}^r[a_k]_{\lambda}}{\prod_{l=1}^s[b_l]_{\lambda}}\frac{C_{\lambda}(\mathbf{x})}{|\lambda|!}\nonumber\\
	=&W^{(s+1)}_{0}(\tilde{N}\beta,\mathbf{b;x}){_rF_s}(\mathbf{a;b;x}).
\end{align}
This completes the proof of Theorem \ref{W0W+-}.
\end{proof}

\begin{rem}
The hypergeometric equation of Gauss-type series $_2F_1(a_1,a_2;b;\x)$ were firstly given by \cite{Yan92} and \cite{Kaneko93} independently. Recently, we proved Corollary \ref{Wxycons} by using (\ref{WO}) and some basic equation system for $_1\mathcal{F}_0(N\beta;\mathbf{x;y})$ in \cite{Liu25}. In Chen and Sahi's original proof of Theorem B and Theorem C \cite{Chen25}, the degree $0$ operators in differential equations were given by the Debiard-Sekiguchi operators but not the symmetric operators $\{\tilde{W}^{(r)}_{0}\}_{r\in\mathbb{N}}$ and $\{W^{(s+1)}_{0}\}_{r\in\mathbb{N}}$ in (\ref{Wop}). Here we explain that these operators are equivalent since they have the same eigenvalues. Meanwhile, this result corrects and proves the conjecture in \cite{Liu25} that the differential operators in the hypergeometric equations can be factored by the Polychronakos' Dunkl operators.
\end{rem}

\subsection{Holonomic system for $_rF_s(\mathbf{a;b;x})$}
By Theorem \ref{Woprep}, the Theorem \ref{W0W+-} means that
\begin{align}
	\sum_{i=1}^Nx_i^{\alpha}\left[\tilde{\W}_{0,i}^{(r)}(\mathbf{a;x})-\W_{-1,i}^{(s+1)}(\tilde{N}\beta,\mathbf{b;x})\right]{_rF_s}(\mathbf{a;b;x})&=0, \quad \alpha=0,1.
\end{align}	
It just is a direct corollary of the following result.
\begin{theo}\label{Hsystem}
The hypergeometric series $_rF_s(\mathbf{a;b;x})$ is the unique solution of the system
\begin{align}\label{hsystem}
	\left[\tilde{\W}_{0,i}^{(r)}(\mathbf{a;x})-\W_{-1,i}^{(s+1)}(\tilde{N}\beta,\mathbf{b;x})\right]F(\mathbf{x})=0, \quad \forall i=1,\cdots,N.
\end{align}
\end{theo}
\begin{proof}
It follows form Proposition \ref{FKFrs} and Lemma \ref{xynonK} that
\begin{align}
	(x_i^{-1}\pi_i-1){_0F_0}(\mathbf{a;b;x}) =0.
\end{align}
By (\ref{OxxO}), (\ref{xOOx}) and (\ref{KrsK00}), we further obtain
\begin{align}
	&\OO^{(s)}(\mathbf{b;x})\left[x_i^{-1}\pi_i\prod_{l=1}^{s}(\pi_i-\tilde{N}\beta+b_{l})\right]{_rF_s}(\mathbf{a;b;x})\nonumber\\
	=&x_i^{-1}\pi_i\OO^{(s)}(\mathbf{b;x}){_rF_s}(\mathbf{a;b;x})\nonumber\\
	=&x_i^{-1}\pi_i\OO^{(r)}(\mathbf{a;x}){_0F_0}(\x)\nonumber\\
	=&\OO^{(r)}(\mathbf{a;x})\left[\prod_{l=1}^{r}(\tilde{\pi}_i+a_{l})x_i^{-1}\pi_i\right]{_0F_0}(\x)\nonumber\\
	=&\left[\prod_{l=1}^{r}(\tilde{\pi}_i+a_{l})\right]\OO^{(r)}(\mathbf{a;x}){_0F_0}(\x)\nonumber\\
	=&\OO^{(s)}(\mathbf{b;x})\left[\prod_{l=1}^{r}(\tilde{\pi}_i+a_{l})\right]{_rF_s}(\mathbf{a;b;x}).
\end{align}
Since $\OO^{(s)}(\mathbf{b;x})$ is invertible, the series ${_rF_s}(\mathbf{a;b;x})$ indeed is the solution of the system (\ref{hsystem}).
\end{proof}

By applying Theorem \ref{Hsystem} recursively and summing over index $i=1,\cdots,N$, we arrive at
\begin{coro}
The hypergeometric series $_rF_s(\mathbf{a;b;x})$ satisfy the following equations
\begin{subequations}
\begin{align}
	\left[\tilde{W}_0^{(kr)}(\mathbf{a}^k;\x)-W_{-k}^{(s+1)}(\tilde{N}\beta,\mathbf{b;x})\right]F(\mathbf{x})&=0,\\
	\left[W_{k}^{(r)}(\mathbf{a};\x)-W_{0}^{(sk+k)}((\tilde{N}\beta,\mathbf{b})^k;\x)\right]F(\mathbf{x})&=0,
\end{align}
\end{subequations}
for any $k\in\mathbb{Z}_{\ge0}$, where $\mathbf{a}^k$ is made up of $k$ parameter tuple $\mathbf{a}$. Especially, the result reduces to Theorem \ref{W0W+-} at $k=1$.
\end{coro}

Following Kaneno's proof of Theorem 4.12 in \cite{Kaneko96}, we give the uniqueness property.
\begin{theo}
The hypergeometric series $_rF_s(\mathbf{a;b;x})$ is the unique solution of the holonomic system (\ref{hsystem}) subject to the following conditions\\
(a) $F(\x)$ is a symmetric function in $x_1,\cdots,x_N$,\\
(b) $F(\x)$ is analytic at the origin with $F(0^N)=1$.
\end{theo}
\begin{proof}
Assuming $S(\x)$ is a  symmetric function solution for the system (\ref{hsystem}) and satisfies $S(0^N)=0$. Let $\x_n=(x_1,\cdots,x_n)$ for $n=1,\cdots,N$. By taking $x_k=0$ for $k=n+1,\cdots,N$, we have that $S(\x_n)=S(x_1,\cdots,x_n,0,\cdots,0)$ also is the solution of the system
\begin{align}
	\mathcal{L}_i(\x_n)F(\x_n)=0, \quad \forall i=1,\cdots,n,
\end{align}
where we denote $\mathcal{L}_i(\x_n)=\tilde{W}_{0,i}^{(r)}(\mathbf{a},\x_n)-\W_{-1,i}^{(s+1)}(\tilde{N}\beta,\mathbf{b};\x_n)$.
Since $S(\x_n)$ is a symmetric function for the variable $\x_n$, we make a formal expansion by Jack polynomial
\begin{align}
	S(\x_n)=\sum_{\lambda\in\mathcal{P}_n}\Omega_{\lambda}\frac{C_{\lambda}(\x_n)}{|\lambda|!}.
\end{align}
By applying the Pieri formula (\ref{Pieri}), the summation constraints $\sum_{i=1}^n\mathcal{L}_i(\x_n)S(\x_n)=0$ gives
\begin{align}
	&\sum_{\lambda\in\mathcal{P}_n}\sum_{i= 1}^{\ell(\lambda)+1}\prod_{l=1}^r(b_l+c_i(\lambda))\Omega_{\lambda^{(i)}}\binom{\lambda^{(i)}}{\lambda}_{\beta}\frac{C_{\lambda^{(i)}}(1^n)}{C_{\lambda}(1^n)}\frac{C_{\lambda}(\x_n)}{(|\lambda|+1)!}\nonumber\\
	=&\sum_{\lambda\in\mathcal{P}_n}\tilde{G}^{(r)}_\lambda(\mathbf{a})\Omega_{\lambda}\frac{C_{\lambda}(\x_n)}{|\lambda|!}.
\end{align}
It further gives the recurrence relation of $\Omega_{\lambda}$ for $\lambda\in\mathcal{P}_n$ as follows
\begin{align}\label{gagai}
	\frac{1}{|\lambda|+1}\sum_{i= 1}^{\ell(\lambda)+1} \prod_{l=1}^r(b_l+c_i(\lambda))\binom{\lambda^{(i)}}{\lambda}_{\beta}\frac{C_{\lambda^{(i)}}(1^n)}{C_{\lambda}(1^n)}
	\Omega_{\lambda^{(i)}}-\tilde{G}^{(r)}_\lambda(\mathbf{a})\Omega_{\lambda}=0.
\end{align}
Now, we intend prove $S(\x_n)\equiv0$ by induction on $n$, the case $n=1$ being immediate from $S(0,\cdots,0)=0$ and (\ref{gagai}). Assuming $S(\x_n)\equiv0$ for $n\le N-1$,  i.e. $\Omega_{\lambda}=0$ for $\ell(\lambda)\le N-1$. For a partition $\lambda$ with $\ell(\lambda)=N-1$, the recurrence relation (\ref{gagai})
shows that $\Omega_{\lambda^{(N)}}$ is a linear combination of $\Omega_{\lambda}$ and $\{\Omega_{\lambda^{(i)}}\}_{i\le N-1}$. It means $\Omega_{\lambda^{(N)}}=0$. By a induction on the value of $\lambda_N$, we get that $\Omega_{\lambda}=0$ for any $\lambda\in\mathcal{P}_N$, i.e., $S(\x)\equiv0$. By Theorem \ref{Hsystem}, we get $_rF_s(\mathbf{a;b;x})$ is the unique solution under the conditions (a) and (b).
\end{proof}

\begin{rem}
Theorem \ref{Hsystem} means that
\begin{align}
	\left[\pi_i\prod_{l=1}^{s}(\pi_i-\tilde{N}\beta+b_{l})-\prod_{l=1}^{r}(\pi_i-\tilde{N}\beta+a_{l})x_i\right]{_rF_s}(\mathbf{a;b;x})=0,
\end{align}
for $i=1,\cdots,N$. In the case of $N=1$, $\pi_i$ reduces to $x\frac{\partial}{\partial x}$. Then the system return to the ODE (\ref{ODE}). It shows that the Dunkl operators play an extremely important role in characterizing Macdonald's $\beta$-deformed hypergeometric series. We believe that the Macdonald's hypergeometric function can be extended to two-parameter deformation and root system $R$ by applying suitable Dunkl operators $\{\pi_{\alpha}:\alpha\in R\}$ and ODE systems $\{\mathcal{L}_i\}_{i=1,\cdots,N}$, such as
\begin{align}
	\mathcal{L}_i(\pi_{\alpha})F=0,\quad i=1,\cdots,N.
\end{align}
Moreover, many special functions associated to general root system in Refs.\cite{Heckman87,Koornwinder,Mac03} and their combinatorial properties may be applied in a unified platform of hypergeometric functions.
\end{rem}

\begin{ex}
The Gauss-type series ${_2F_1}(a,b;c;\mathbf{x})$ satisfies the following holonomic system
\begin{align}\label{21system}
	x_i(x_i-1)\partial_i^2F+\left[(a+b+2-\tilde{N}\beta)x_i-(c-\tilde N\beta+1)\right]\partial_i F+abF\nonumber\\
	+\beta \sum_{j\ne i}\left(\frac{x_i-1}{x_i-x_j}x_i\partial_i-
	\frac{x_j-1}{x_i-x_j}x_j\partial_j\right)F=0, \qquad i=1,\cdots, N.
\end{align} 
\end{ex}
In Ref.\cite{Kaneko93}, Kaneko derived the holonomic system (\ref{21system}) for a Selberg-type integral and then proved that the integral is identical to  the hypergeometric series ${_2F_1}(a,b;c;\mathbf{x})$ by uniqueness theorem. Integral representations of these hypergeometric functions have been studied extensively in Refs.\cite{Baker97,Brennecken22,Yan92}, so we do not discuss them further here.

\end{document}